\documentclass[12pt]{amsart}    
\usepackage{amssymb,amsmath,eepic,verbatim, bm, bbm}
\usepackage{amsfonts}
\usepackage{amsbsy}
\usepackage{mathtools}
\usepackage{enumitem}

\usepackage{relsize}

\newtheorem*{thm*}{Theorem}
\newtheorem*{lem*}{Lemma}
\newtheorem*{prop*}{Proposition}
\newtheorem*{cor*}{Corollary}
\newtheorem*{claim*}{Claim}

\newtheorem{thm}{Theorem}
\newtheorem{lem}{Lemma}
\newtheorem{prop}{Proposition}
\newtheorem{cor}{Corollary}
\newtheorem{ex}{Example}
\newtheorem{exs}{Examples}

\usepackage{amsthm}

\theoremstyle{definition}

\numberwithin{thm}{section}
\numberwithin{lem}{section}
\numberwithin{prop}{section}
\numberwithin{cor}{section}
\numberwithin{ex}{section}
\numberwithin{claim}{section}

\newcommand{\define}{\emph}

\newcommand{\C}{\mathbb{C}}

\newcommand{\Z}{\mathbb{Z}}
\newcommand{\Q}{\mathbb{Q}}

\newcommand{\on}{\operatorname}
\newcommand{\Gr}{ \on{Gr} }
\newcommand{\Fl}{ \on{Fl} }

\newcommand{\bGr}{\bm{Gr}}

\newcommand{\bull}{ {\scriptscriptstyle{\bullet}}  }

\renewcommand{\min}{{\text{min}}}

\newcommand{\Pf}{ \on{Pf} }

\newcommand{\rank}{ \on{rank} }
\newcommand{\Hom}{ \on{Hom} }
\newcommand{\Det}{ \on{Det} }
\newcommand{\Ker}{ \on{Ker} }

\newcommand{\id}{ \on{id} }

\newcommand{\del}{ \partial }
\newcommand{\nearr}{ \mathsmaller{\nearrow} }

\let\emph\relax 
\DeclareTextFontCommand{\emph}{\bfseries\em}

\newcommand{\enS}{\mathsf{S}}              
\newcommand{\tenS}{\mathbf{S}}
\newcommand{\tenG}{\mathbf{G}}              
\newcommand{\Sgp}{\mathrm{Perm}}   
\newcommand{\LLambda}{\Lambda}    
\newcommand{\GGamma}{\Gamma}      
\newcommand{\bLambda}{\bm{\Lambda}}    
\newcommand{\bGamma}{\bm{\Gamma}}   
\newcommand{\tomega}{\bm{\omega}}    
\newcommand{\balpha}{\bm\alpha}          

\begin{document}
\title{Schubert Polynomials in Types A and C}
\date{\today}
\author{David Anderson}
\email{anderson.2804@osu.edu}
\address{Department of Mathematics,
The Ohio State University,
Columbus, Ohio, 43210}
\thanks{DA is partially supported by NSF CAREER DMS-1945212.}

\author{William Fulton}
\email{wfulton@umich.edu}
\address{Department of Mathematics,
University of Michigan,
Ann Arbor, Michigan  48109}

\maketitle

\begin{abstract}
Enriched versions of type A Schubert polynomials are constructed with coefficients in a polynomial ring in variables $c_1, c_2, \ldots$.  Specializing these variables to $0$ recovers the double Schubert polynomials of Lascoux and Sch\"utzenberger; specializing them to certain power series recovers the back-stable double Schubert polynomials of Lam, Lee, and Shimozono; specializing them to Schur Q-polynomials relates them to the type C double Schubert polynomials of Ikeda, Mihalcea, and Naruse.  Many formulas for classical Schubert polynomials generalize to this setting.  They give, and are characterized by, formulas for degeneracy loci.
\end{abstract}

\tableofcontents

\section*{Introduction}

Lascoux and Sch\"utzenberger constructed  single and double Schubert polynomials  in polynomial rings  $\Z[x_+]$ and $\Z[x_+,y_+]$  in variables $x_i$ and $y_j$  for positive $i$ and $j$, with integer coefficients.   Billey and 
Haiman \cite{BH} realized that to generalize to type C
required taking coefficients in another ring $\GGamma$, generated by certain Schur Q-polynomials $q_i$.  
Their Schubert polynomials were generalized to double polynomials by Ikeda, Mihalcea, and Naruse 
\cite{IMN}, cf.~\cite{AF1} and \cite{BKT}. (For types B and D a similar ring $\GGamma'$, generated by elements $p_i$, was used, 
but we will not discuss these here.)

Having coefficients in $\GGamma$  makes the type C  story richer than the type A story. 
For example, in type C
the Stanley symmetric functions can be obtained by setting the polynomial variables 
equal to zero, while this fails in type A; and the ring $\GGamma$ controls intersection rings of isotropic Grassmannian bundles.  Our main object in 
this paper is to construct a version of  \emph{enriched Schubert polynomials} $\enS_w$ in type 
A that are in a ring $\Lambda[x,y]$.  The ring $\Lambda$ 
is, in fact, a polynomial ring  $\Z[c_1,c_2,\ldots]$, with variables $c_k$ of degree $k$. The 
classical double Schubert polynomials are recovered by setting all these $c_k$ to $0$, and the type A Stanley polynomials are recovered by setting all the $x$'s and $y$'s to $0$.  Most of the classical properties of classical Schubert polynomials, e.g.~in \cite{Mac2}, now have enriched analogues in all four types.

It has been known at least since \cite{BH} that the type A
Stanley symmetric polynomial for $w \in \Sgp_{\infty}$ maps to the type C Stanley symmetric polynomial for the 
same $w$ regarded as a signed permutation.  We will see that the homomorphism from $\Lambda$ to $\GGamma$ that maps  $c_i$ to $q_i$ sends 
the type A Schubert polynomial for a usual permutation to the type C polynomial for that permutation.\footnote{This is implicit in \cite{LLS2}, but our approach provides a clear geometric reason, as well as a stronger result for the ``twisted'' polynomials described below.}

The Stanley symmetric polynomials are motivated by comparing classical Schubert polynomials for a permutation $w$ with those of the permutation 
$\gamma^m(w) = 1^m \times w$; that is, the permutation taking $i$ to $i$ for $i\leq m$ and $m+i$ to $m+w(i)$ for $i > 0$, cf.~\cite{BJS}.   The polynomials constructed here also depend on this ``back-stable" property.  
In fact, these polynomials are very close to the back-stable Schubert polynomials defined and developed in~ \cite{LLS1}.  These back-stable Schubert polynomials are in a ring  
$\overleftarrow{R}[x,a]$, 
which is generated by certain infinite series $p_k(x||z)$, as well as $x$ and $a$ variables.  
There is a homomorphism from $\Lambda$ to  
$\overleftarrow{R}$, which, with $y_i$ mapping to $-a_i$, sends the Schubert polynomials  $\enS_w$ to their back-stable 
Schubert polynomials $\overleftarrow{\mathfrak{S}}_w(x;a)$.  Since this homomorphism is an embedding, our polynomials are determined by theirs. 
In spite of this, we believe the simplicity of our approach has some advantages.

Among the many important innovations in \cite{LLS1} is the generalization of Schubert polynomials from permutations of a finite set of positive integers, to all permutations in $\Sgp_\Z$ that permute any finite set of 
integers.  That generalization is included here, producing a Schubert polynomial $\enS_w$ in $\Lambda[x,y]$ for any 
$w$ in $\Sgp_\Z$, where now $x$ and $y$ stand for variables $x_i$ and $y_i$ for all integers $i$.

We give \textit{determinantal formulas} for $\enS_w$ for all vexillary $w$, which includes an expression for them in terms of the basis of $\Lambda$ by \textit{Schur polynomials}
\[S_\lambda(c) = \Det(c_{\lambda_i+j-i}).\]  
Unlike in the classical case, these $\enS_w$ are always irreducible polynomials, for any vexillary $w$.  All the Schubert polynomials are determined by the vexillary ones, using difference operators (or transition formulas) relating them.  Identities from \cite{AF2} give explicit formulas for these difference operators and translation operators, as well as the difference operators in type C.

There is a \textit{twisted} version of this story, with a parameter $z$ that corresponds to the first Chern class of a line bundle.  This generalization is straightforward in type A, with the ring $\Lambda$ simply replaced by $\bLambda = \Lambda[z]$; we denote these twisted enriched Schubert polynomials by $\tenS_w$.  The twists in the other classical types are less trivial, and the fact that the twisted polynomials in types A and C are compatible is algebraically more surprising.

As with the classical double Schubert polynomials \cite{F2}, the Schubert polynomials defined here give formulas for degeneracy loci, and they are determined by these formulas. 
In type C one has a vector bundle $V$ with a symplectic form $V \times V \to L$ for some line bundle $L$, and one has maximal isotropic subbundles $E$ and $F$, each in the middle of complete isotropic filtrations.  The ring $\bGamma$ used for the (twisted) type C story is  a residue ring $\bLambda/I$, for a certain ideal $I$.  The $c_k$  in type C map to the Chern classes  $ c_k(V - E - F) $.

The natural corresponding type A setting is to have two vector bundles $V$ and $W$, with a bilinear mapping $V \times W \to L$.  Assume we have filtrations 
\[ 
 \cdots  E_1 \subset E = E_0 \subset E_{-1} \subset \cdots \subset V, \;\;  \cdots  F_1 \subset F = F_0 \subset F_{-1} \subset \cdots \subset W ,
\]
with $E$ and $F$ of the same rank; the ranks of  $E_i$ and $F_i$ are equal to $\rank(E) - i$ for all $i$.  For $w \in \Sgp_\Z$ define a degeneracy locus $\Omega_w(E_\bull,F_\bull)$ by requiring 
\[
  \dim( \Ker (F_q \to  \Hom(E_p,L) )) \geq \, \#\{a \leq p \mid w(a) > q \} 
\]
for all $p$ and $q$.\footnote{The left side could be called the \textit{nullity} of the bilinear map, following the classical notation that the nullity of an $m$ by $n$ matrix is the dimension of the kernel of the corresponding linear map from $n$-space to $m$-space.}
The enriched Schubert polynomials are characterized by the fact that, whenever $\Omega_w(E_\bull,F_\bull)$  has the expected codimension $\ell(w)$, its class is given by the formula $\tenS_w(c,x,y)$, where now 
\[
 c_k \mapsto c_k(E^*\otimes L - F), \; x_i \mapsto c_1(E_{i-1}/E_i), \; y_j \mapsto c_1(F_{j-1}/F_j).
\]
They are uniquely determined by the special case when the bundles are all direct sums of line bundles, in which case these can be regarded as variations on the matrix determinantal loci studied in~\cite{KM}.  Many of their properties become evident from this geometric interpretation.  Other cases of this degeneracy locus formula were observed by Knutson and studied by Pawlowski \cite{P}.

The main text of this paper states properties of these enriched Schubert polynomials. We expect that readers with varying backgrounds will readily supply proofs from those backgrounds.  Many of the properties follow routinely from the corresponding properties of classical Schubert polynomials.  Many also follow from the formulas for degeneracy loci.  And some can be deduced from corresponding results of~\cite{LLS1}. For brevity and to order the results naturally, most of these are omitted here; systematic details will be found in \cite{AF4}.  Some proofs which involve new ideas are sketched in a final section. 

There is another relation to geometry, also considered in \cite{LLS1}, for which these enriched Schubert polynomials are equivariant classes of Schubert varieties in infinite flag varieties.  This also gives proofs of many of the results here.  Since this involves subtleties of limits and colimits of varieties, it is left to a separate paper \cite{A}.

All our formulas involve only polynomials with a finite number of variables, and with integer coefficients.

\medskip
\noindent\emph{Question.} Buch, in \cite{Bu}, constructed a filtered ring, which we denote by $A$ to 
avoid conflicting with other notation here, whose associated graded ring is 
the ring $\LLambda$.  It is natural to ask if there are \textit{twisted enriched 
Grothendieck polynomials} $\tenG_w$ in $A[z][x,y]$, which represent the 
classes of the degeneracy loci $\Omega_w$ in the Grothendieck ring of 
vector bundles, and whose leading terms are the enriched Schubert 
polynomials $\tenS_w$ of this paper.  The possibility of $K$-theoretic analogues of the 
back-stable Schubert polynomials was announced in \cite{LLS1}.

\medskip
\noindent\emph{Acknowledgements.} This paper is intended as an introduction to, and variations on -- as well as a tribute to -- the groundbreaking work of Lam, Lee, and Shimozono~\cite{LLS1}.  We should also again thank Kazarian~\cite{K}, who saw that using bundles of ranks near the middle gave the geometry behind type C, B, and D Schubert polynomials; this paper is the result of carrying out a similar idea in type A.

\medskip
\noindent\emph{Remark.}  The point of view in this paper comes from geometry: the coefficient rings $\LLambda$ and $\GGamma$ for Schubert polynomials are regarded as \textit{source rings}, whose generators are free to be mapped to Chern classes of (virtual) vector bundles, or to be otherwise specialized.  The opposite view of them as \textit{target rings} is common in combinatorics, which often identifies them with, or embeds them in, rings of symmetric functions.  For us $\LLambda$ is the polynomial ring $\Z[c_1, c_2, \ldots]$ in variables $c_i$.  Of course, this ring can be identified with the ring of symmetric functions, with  $c_i$ mapping to the corresponding complete homogeneous symmetric function, and $S_\lambda(c)$ defined here becomes the Schur function $s_\lambda$.  The variables used for the symmetric functions will not appear here.

\medskip
\noindent\emph{Notation.}  Geometry also motivates some of our notation.  In geometry one writes $c(E) = \sum c_i(E)$  for the total Chern class of a vector bundle $E$; this notation has replaced an earlier cumbersome notation which used Chern polynomials  $c_t(E) = \sum c_i(E) t^i$.  One writes now  $c(E - F) = c(E)/c(F)$ in place of the corresponding division of power series.  Given elements $c_1, c_2,\ldots $ in a commutative ring, we will often write $c = \sum c_i$ as an abbreviation for the series $\sum c_i t^i$, with the standard convention that  $c_0 = 1$, and call this a \emph{Chern series}.  We use the corresponding group notation, so $c = a \cdot b$  means that $c_k = \sum_{i+j = k} a_i b_j$.  We write similarly $c = \prod_{i=1}^m \frac{1+x_i}{1+y_i}$ in place of $\sum c_k t^k = \prod_{i=1}^m \frac{1+x_it}{1+y_it}$.  (Letters $a, b, c$ at the beginning of the alphabet will be used for these Chern series, and letters like  $x_i, y_j, z$ at the end for variables of degree $1$ that appear in such equations.)

For any Chern series $c(1), c(2), \ldots, c(n)$ and partition $\lambda$ of length at most $n$,  define the \emph{Schur determinant} 
 \[S_\lambda(c(1), \ldots, c(n)) \, = \, \Det(c(k)_{\lambda_k+l-k})_{1 \leq k,l \leq n}.\]
 For a pair of partitions $\mu \subset \lambda$, set
  \[S_{\lambda/\mu}(c(1), \ldots, c(n)) \, = \, \Det(c(k)_{\lambda_k-\mu_l+l-k})_{1 \leq k,l \leq n}.\]

 For any commutative ring $R$, $R[x,y]$ denotes the polynomial ring in variables $x_i$ and $y_j$  for all integers $i$ and $j$, and  $R[x_+,y_+]$ is the polynomial ring in variables $x_i$ and $y_j$ for $i$ and $j$ positive integers.
 
The classical single and double Schubert polynomials of Lascoux and 
Sch\"utzenberger, for $w$ in $\Sgp_{+}$ are denoted by $\mathfrak{S}_w(x)$ and 
 $\mathfrak{S}_w(x;y)$; these are in polynomial rings  $\Z[x_+]$ 
 and $\Z[x_+,y_+]$.

We follow  \cite{LLS1} for symmetric group notation, except that we write $\Sgp$ instead $S$ for permutation groups to avoid a plethora of $S$'s.
The group of permutations of the integers that fix all but a finite 
number of integers is denoted $\Sgp_\Z$.  It is generated by the 
elements  $s_i$ that transpose $i$ and $i+1$, for $i \in \Z$.  
It contains the subgroup 
$\Sgp_{+} = \Sgp_\infty = \bigcup \Sgp_n$ that permutes the positive integers, 
the subgroup $\Sgp_{-}$ that permutes the negative integers, as 
well as the subgroup $\Sgp_{\neq 0} = \Sgp_{-} \times \Sgp_{+}$.  

Lengths and reduced words are defined for $w$ in $\Sgp_\Z$  just as for 
$w$ in $\Sgp_{+}$; the length or inversion number of $w$ is denoted $\ell(w)$.  
The notation $u \cdot v \stackrel{.}{=} w$ means the product of 
$u$ and $v$ is $w$ in $\Sgp_\Z$, with the condition that $\ell(v)+\ell(v) = \ell(w)$.

We use a \emph{one-line notation} for $w \in \Sgp_\Z$: select a sequence of $n$ consecutive integers containing all integers not fixed by $w$, and  write them in the order that $w$ puts them.  For example, $w = 3 \, -\!1 \,\, 2 \,\, 0 \, -\!2 \,\, 4 \, 1$  permutes the integers from $-2$ to $3$, mapping $-2$ to $3$, $-1$ to $-1$, $0$ to $2$, $1$ to $0$, $2$ to $-2$, $3$ to $4$, and $4$ to $1$, and fixing all other integers.  Equivalently, if $w$ fixes all integers outside an interval $[p+1,p+n]$, we write $w = w_1  \cdots w_n$, where  $w_i = w(p+i)$ for $1 \leq i \leq n$, and $p+1 = \min(w_1, \ldots, w_n)$.  
This representation  is unique if $w_1$ is not the 
minimum and $w_n$ is not the maximum in these sets of integers; 
if either of these holds, such integers can be omitted.  Note that 
for $w$ in $\Sgp_{+}$, this agrees with the usual one-line notation, except that 
the usual notation always uses an interval $[1,n]$.  So here one could 
write $w = 4\,3$ instead of $1\,2\,4\,3$ for  $w = s_3$.

There is an automorphism $\gamma$ of $\Sgp_\Z$ that maps $s_i$ to 
$s_{i+1}$.  The automorphism  $\gamma^m$, for any integer 
$m$, maps $s_i$ to $s_{i+m}$.  For $w$ in $\Sgp_{+}$ and $m$ positive,  
$\gamma^m(w)$ is the permutation often denoted $1^m \times w$.  
For any integer $m$, the one-line notation for $\gamma^n (w)$  is obtained from that of $w$  by adding  $m$  to each entry.  

There is an involution $\omega$ of $\Sgp_\Z$ that takes $s_i$ to $s_{-i}$.  In one-line notation, if $w= w_1\;\cdots\; w_n$, then $\omega(w) = 1-w_n \; \cdots \; 1-w_1$.

\section{Construction of the Schubert Polynomials}

Let $\LLambda = \Z[c]= \Z[c_1, c_2, \ldots]$ be the polynomial ring 
in variables $c_k$, for positive integers $k$; this is a graded ring, with 
$c_k$ having degree $k$.  
The polynomial ring $\LLambda[x,y]$ has automorphisms, 
denoted $\gamma^m$,  for 
each integer $m$, with $\gamma^0$  the identity automorphism, and 
\[\gamma^m \circ \gamma^n = \gamma^{m+n}\] 
for all $m$ and $n$. 
First, $\gamma^m(x_i) = x_{m+i}$  and $\gamma^m(y_i) = y_{m+i}$  
for all integers $m$ and $i$.  The action on the generators $c_k$ of 
$\LLambda$ is defined as follows.
  For $m = 0$, $\gamma^m(c_k) = c_k$.  For  $m>0$, 
\[ \gamma^m(c_k) = 
\sum_{p+q+r=k} c_p \, h_q(x_1, \ldots, x_m)\, e_r(y_1, \ldots, y_r),\]
where $h_q$ and $e_r$ denote the complete and elementary symmetric 
in the indicated variables.  For $m < 0$,
\[\gamma^m(c_k) = 
\sum_{p+q+r=k} c_p \, e_q(-x_{m+1}, \ldots, -x_0)\, h_r(-y_{m+1}, \ldots, -y_0).\]
Using the multiplicative notation, this can be written

\[
\gamma^m(c) = c \cdot  \prod_{i=1}^m \frac{1+y_i}{1-x_i} \cdot \prod_{i=m+1}^0 \frac{1-x_i}{1+y_i} ,
\]
with the understanding that only one of these products is used, depending on whether $m$ is positive or negative.

The following result can be used for explicit calculations:
\begin{prop}\label{gamma}
For $m > 0$, 
\[
\gamma^m(S_\lambda(c)) = \sum_{\mu \subset \lambda} \sum_{T} (x,y)^T S_\mu(c), 
\]
the second sum over all tableaux $T$ on the skew shape $\lambda/\mu$  with entries 
$1' < 1 < 2' < 2 < \dots < m' < m$, weakly increasing along rows and down columns, with no $k$ repeated in a column and no $k'$ repeated in a row; and $(x,y)^T = \prod x_k^{\# \{k \in T\}} \prod y_k^{\#\{ k' \in T\}}$.
\end{prop}

\begin{thm}\label{DefThm}
For every $w$ in $\Sgp_\Z$, there is a unique polynomial  $\enS_w = \enS_w(c,x,y)$ in $\LLambda[x,y]$
such that for all positive integers $m$ such that $\gamma^m(w)$ is in $\Sgp_{+}$, 
\[
 \gamma^m(\enS_w)(0,x,y) = \mathfrak{S}_{\gamma^m(w)}(x,-y) .
\]
\end{thm}

In fact, if $w$ is in $\Sgp_{+}$, $\enS_w$ is determined by this specialization for any one $m$ that is sufficiently large: any $m \geq \ell(w)/2$ will do.  If $w$ is not in  $\Sgp_+$, and $n$ is the smallest integer not fixed by $w$, it suffices to take  $m \geq \ell(w)/2  - n$.

Each $\enS_w$ is a homogeneous polynomial of degree $\ell(w)$, where $c_k$ has degree $k$ and each $x_i$ and $y_j$ has degree $1$.  If $w(i) = i$ for all $i$ not in an interval $[-m,n]$, for $m \geq 0$, $n > 0$, the variables $x_i$ and $y_i$  only appear if $i$ is in $(-m,n)$. 

\begin{exs} Set $\alpha_{i j} = x_i + y_j$. 
\begin{enumerate}
\item For $w$ the identity permutation, $\enS_w = 1$.
\item For $k = 0$, $\enS_{s_k} =  c_1$; for $k > 0$, $\enS_{s_k} =  c_1 + \sum_{i=1}^k\alpha_{i i}$; for $k < 0$, 
$\enS_{s_k} =  c_1 - \sum_{i=k+1}^0\alpha_{i i}$.
%
%
%
%
\item   $\enS_{2\,1\,4\,3}  = c_1^2+(2\alpha_{1 1}+\alpha_{2 2} + \alpha_{3 3})c_1 +
\alpha_{1 1}(\alpha_{1 1}+\alpha_{2 2} + \alpha_{3 3})$.
%
%
\end{enumerate}
\end{exs}

Since the polynomials $S_\lambda(c) = \Det(c_{\lambda_i + j - i})$ form a basis for $\LLambda$, these Schubert polynomials have unique expressions in terms of this basis.  The table at the end of the paper includes all of these for all $w$ permuting $0$, $1$, $2$, and $3$.

These enriched Schubert polynomials specialize to classical double Schubert polynomials: for $w$ in $\Sgp_{+}$, 
\[
\enS_w(0,x,y) = \mathfrak{S}_w(x,-y).
\]

They are \emph{stable} in the usual sense: they are defined for any 
permutation $w$ of an interval $[-m,n]$ of integers, with $m \geq 0$ and $n > 0$, as a polynomial in $\LLambda[x_{-m+1},\ldots, x_{n-1},y_{-m+1},\ldots, y_{n-1}]$, and the result is independent of the choice of $m$ and $n$.  

They satisfy a \emph{duality} or \emph{inverse} property.  Define variables $\omega(c)_k$ to be the Schur polynomials: $\omega(c)_k = S_{1^k}(c)$ for the partition consisting of $k$ $1's$; equivalently, $\omega(c) = 1/(1-c_1+c_2 -\cdots)$.  This is the well-known involution of $\LLambda$ that takes any $S_\lambda(c)$ to $S_{\lambda'}(c)$, where $\lambda'$ is the conjugate of $\lambda$, cf.~\cite{Mac1}.

\begin{prop}\label{duality}  For any $w$ in $\Sgp_\Z$,  
\[
  \enS_{w^{-1}}(c,x,y) = \enS_w(\omega(c),y,x). 
\]
\end{prop}

There is an involution $\omega$ of $\LLambda[x,y]$ that takes $x_i$ to $-x_{1-i}$, $y_i$ to $-y_{1-i}$, and $\omega(c)$ as in the preceding proposition.  The following is equivalent to Proposition~4.17 of \cite{LLS1}.

\begin{prop}\label{omega} For any $w$ in $\Sgp_\Z$,
\[\enS_{\omega(w)} = \omega(\enS_w).\]
\end{prop}

\section{Back-Stability}

These polynomials are  \emph{back-stable}, in the following sense (cf.~\cite[Cor.~4.4]{LLS1}):
\begin{prop}\label{p.back-stable}
For any integer $m$ and any $w$ in $\Sgp_\Z$,  
\[\enS_{\gamma^m(w)} \, = \, \gamma^m(\enS_w) .\]
\end{prop}

Lam, Lee, and Shimozono \cite{LLS1} construct a ring $\LLambda(x||a)$, which is an algebra over  $\Q[a]$, the polynomial ring over $\Q$ in variables $a_i$, $i \in \Z$.  In fact, it is a polynomial ring over $\Q[a]$, with generators  $p_k(x||a)= \sum_{i \leq 0} \, x_i^k - a_i^k$, power series in two sets of variables.  Their back-stable Schubert polynomials $\overleftarrow{\mathfrak{S}}_w(x;a)$ are in the algebra 
\[ \overleftarrow{R}(x;a)  = \LLambda(x||a) \otimes_{\Q[a]} \Q[x,a] ,\]
where $\Q[x,a]$ is the polynomial ring in variables $x_i$ and $a_i$ for all integers $i$.
For any partition $\lambda$, set
\[p_\lambda(x||a) = \prod_{j=1}^{\ell(\lambda)}p_{\lambda_j}(x||a) = \prod_{i>0} p_i(x||a)^{m_i(\lambda)}, \] 
where $m_i(\lambda)$ is the multiplicity with which $i$ occurs in $\lambda$.  Let $z_\lambda = \prod_{i>0} i^{m_i(\lambda)}\cdot m_i(\lambda)!$.  Define a homomorphism $\rho \colon \LLambda \to \LLambda(x||a)$ by defining $\rho(c_k)$ by the formula
\[\rho(c_k) =  \sum_{|\lambda| = k} \,\,\frac{1}{z_\lambda}\,  p_\lambda(x||a).\]
This extends to a homomorphism  $\rho \colon \LLambda[x,y] \to \overleftarrow{R}(x;a) $ that takes $x_i$ to $x_i$ and $y_i$ to $-a_i$.

\begin{prop}\label{LLS1}  For any $w$ in $\Sgp_\Z$, $\rho(\enS_w) = \overleftarrow{\mathfrak{S}}_w(x;a)$.
\end{prop}

\section{Difference Operators}

There are \emph{difference operators} $\del_i = \del_i^x$ which are endomorphisms 
of $\LLambda[x,y]$, linear over $\Z[y]$.  There are corresponding algebra 
automorphisms $s_i$, with 
\[ \del_i(f) = \frac{f - s_i(f)}{x_i - x_{i+1}},\]
and they satisfy the Leibniz rule
$ \del_i(f \cdot g) = \del_i(f) \cdot g + s_i(f) \cdot \del_i(g)$.
They act as usual on the $x$ variables, with $s_i$ interchanging $x_i$ and $x_{i+1}$,
so $\del_i(x_i) = 1$ and $\del_i(x_{i+1}) = -1$, and $\del_i(x_j) = 0$ for all other 
integers $j$. For $i \neq 0$, the $c_k$ act as scalars; but 
$ \del_0(c_k) =  \sum_{i+j =k-1} x_1^i c_j $.

\begin{prop}\label{del0}
\[\del_0(S_\lambda(c)) \, = 
\, \sum (x_1)^{v(\lambda/\mu)} \, (-x_0)^{h(\lambda/\mu)} \, (x_1 - x_0)^{k(\lambda/\mu)-1} S_\mu(c).\]
the sum over all $\mu$ obtained from $\lambda$ by removing a nonempty border strip (from its Young diagram), with $v(\lambda/\mu)$ the number of vertical lines between border 
boxes,  $h(\lambda/\mu)$ the number of horizontal lines 
between border boxes,  and $k(\lambda/\mu)$ the number of connected components in the border strip.
\end{prop}

For example, writing $S_\lambda$ for $S_\lambda(c)$,
\begin{align*}
\del_0( S_{(4,2)} ) &= S_{(4,1)} + S_{(3,2)} + x_1\,S_{(4)} + x_1\,S_{(2,2)} + (x_1-x_0)\,S_{(3,1)}\\
& \qquad + x_1(x_1-x_0)\, S_{(3)} +  x_1(x_1-x_0)\, S_{(2,1)} + x_1^2(x_1-x_0)\, S_{(2)} \\
& \qquad + (-x_0)x_1^2\, S_{(1,1)} + (-x_0) x_1^3 \, S_{(1)}.
\end{align*}

There are corresponding difference operators  $\del_i^y$  defined by
interchanging the roles of the $x$ and $y$ variables.

The difference operators are compatible with the translation operators $\gamma^m$:
\[
\gamma^m \circ \del_i = \del_{i+m} \circ \gamma^m
\]
for all integers $m$ and $i$, where  $\del_i$ denotes either $\del_i^x$ or $\del_i^y$.

\begin{thm}\label{DiffOp}  The polynomials $\enS_w$, for $w$ in $\Sgp_\Z$, satisfy and are
 uniquely determined by the following properties:
 \begin{enumerate} 
 \item $\enS_w$ is homogeneous of degree equal to the length $\ell(w)$ of $w$. 
 \item $\enS_{\id} = 1$.
\item For any $w$ and integer $i$, 
\[ \del_i^x(\enS_w) \, = \begin{cases} \enS_{w s_i}  & \text{if} \;\; \ell(w s_i) < \ell(w) \\
0 & \text{if} \;\; \ell(w s_i) > \ell(w) ;
\end{cases} \]
and
\[ \del_i^y(\enS_w) \, = \begin{cases} \enS_{s_i w}  & \text{if} \;\; \ell(s_i w) < \ell(w) \\
0 & \text{if} \;\; \ell(s_i w) > \ell(w) .
\end{cases} \]
\end{enumerate}
\end{thm}

\noindent
(See also \cite[Thm.~4.6]{LLS1}, where the last condition is replaced by the requirement on specializations of the polynomials.)

\section{Determinantal Formula for Vexillary Permutations}

A permutation  $w$  in  $\Sgp_\Z$  is \emph{vexillary} if it avoids the pattern 
$2\, 1\, 4 \,3$.  This notion is preserved by the translation operators $\gamma$, so each is 
a translate of an ordinary vexillary permutation in $\Sgp_{+}$.  We describe a way to construct these explicitly, and use the 
ingredients to give a determinantal formula for their Schubert polynomials. 
Each is associated with a \emph{type A triple} $\tau = (k_\bull,p_\bull,q_\bull)$ of sequences, each of the same 
length $s$, called the \emph{length} of $\tau$.  The first sequence $k_\bull$
consists of $s$  positive integers in strictly increasing order:
\[
0 <  k_1 < k_2 < \cdots  < k_s .
\]  
The sequences $p_\bull$ and $q_\bull$ are weakly increasing and decreasing
sequences of arbitrary integers:
\[
 p_1 \leq p_2  \leq \cdots  \leq p_s \;\; \text{and} \;\; q_1 \geq q_2 \geq \cdots \geq q_s .
\]
Define integers $l_i$ by the formula $l_i = q_i - p_i + k_i$.  The condition to be a triple is that 
\[
  l_1 > l_2 > \cdots > l_s > 0 .  \tag{*}  
\]
A triple $\tau$ determines a partition $\lambda = \lambda(\tau)$, whose 
Young diagram has corners at $(k_i,l_i)$, $1 \leq i\leq s$.  That is, 
$\lambda$ has length $k_s$ and $s$ distinct parts, and for $1 \leq k \leq k_s$, $\lambda_k = l_i$ for $i$ minimal with  $k_i \geq k$.  

The vexillary permutation $w = w(\tau)$ can be defined as the permutation of minimal 
length such that, for $1 \leq i \leq s$,
\[
 \# \{a \leq p_i \mid w(a) > q_i\} = k_i .
\]
The permutation $w(\tau)$ is in $\Sgp_{+}$ exactly when each $p_i$ and  $q_i$ is 
positive and each $r_i$ is nonnegative, where $r_i = p_i - k_i = l_i - q_i$.  In this case $w$ is of minimal length with 
$\{a \leq p_i \mid w(a) \leq q_i\} = r_i $, as in \cite{F2}.  The identity permutation corresponds to the 
empty $\tau$ of length zero.

For example $\tau = \left( (2,3,5), (1,1,3), (2,0,-1) \right)$ has partition $\lambda = (3^2,2,1^2)$, and permutation $w = 1 \, 3  \, 4\, 0 \, 2  -\!\!1$.

Define the \emph{translate} $\gamma^m(\tau)$ by leaving $k_\bull$ unchanged, but adding 
$m$ to each $p_i$ and $q_i$.  Then $\lambda(\gamma^m(\tau)) = \lambda(\tau)$ and $w(\gamma^m(\tau))= \gamma^m(w(\tau))$. 

There is a \emph{conjugate} triple $\tau'$, with $k'_\bull = (l_s, \ldots, l_1)$, $p'_\bull = (q_s, \ldots, q_1)$, and $q'_\bull = (p_s, \ldots, p_1)$; so $\l'_\bull = (k_s, \ldots, k_1)$.  Then $\lambda(\tau') = \lambda(\tau)'$ and $w(\tau') = w(\tau)^{-1}$.

There is also an \emph{opposite} $\omega(\tau) = (k_\bull, -q_\bull, -p_\bull)$, with $\lambda(\omega(\tau)) = \lambda(\tau)$ and $w(\omega(\tau)) = \omega(w(\tau))$.

Each vexillary  $w = w(\tau)$ has a Schur-type determinantal formula for its 
Schubert polynomial.  For this, we need to define Chern series  $a(p,q)$ with 
coefficients in $\Z[x,y]$ .  Using the multiplicative notation, these are 
defined by 
\[
 a(p,q) \, =  \, \frac{\prod_{i=p+1}^0(1-x_i) \prod_{j=1}^q (1+y_j)}{\prod_{i=1}^p(1-x_i) \prod_{j=q+1}^0 (1+y_j)}. 
\]
Note that at most one of each of the products involving $x$'s or $y$'s can be present; i.e.,~this is shorthand for four formulas, one for each of the possibilities of positive or negative $p$ or $q$.\footnote{Another expression for this is $a(p,q) = \frac{\prod_{i \leq 0}(1-x_i) \prod_{j \leq q} (1+y_j)}{\prod_{i \leq p} (1-x_i) \prod_{j \leq 0} (1+y_j)}$, where the products are over indices greater than some large negative integer, or greater than minus infinity if infinite products are used.}  
Define, for  $1 \leq k \leq k_s$,  $c(k)$ to be the Chern series $c \cdot a(k)$, where $a(k) =   a(p_i, q_i)$, with $i$ minimal so that $k_i \geq k$.  
 
\begin{thm}\label{Vex}
For any $\tau$, with partition $\lambda = \lambda(\tau)$,  the vexillary permutation $w = w(\tau)$ 
has Schubert polynomial 
\[
\enS_w  \, = \, S_{\lambda}(c(1), \ldots, c(k_s)) \, = \, \Det(c(k)_{\lambda_k+l-k}). 
\]
\end{thm}

\begin{cor}\label{vexcor}
For any $\tau$, the vexillary permutation $w = w(\tau)$ 
has Schubert polynomial 
\[
 \enS_{w} = \sum S_{\lambda/\mu}(a(1),\ldots, a(k_s)) S_\mu(c). 
\]

\end{cor}

\begin{cor}\label{irred}
The Schubert polynomial $\enS_w$ of any vexillary $w$ is an irreducible polynomial in the polynomial ring $\LLambda[x,y]$. 
\end{cor}
Note that $\enS_{2\,1\,4\,3} = (c_1+ \alpha_{1 1})(c_1+\alpha_{1 1}+\alpha_{2 2} + \alpha_{3 3})$, where $\alpha_{i i} = x_i + y_i$, so not all enriched Schubert polynomials are
irreducible.

As in \cite{LLS1}, there is a permutation $w_\lambda$ in $\Sgp_\Z$ corresponding 
to each partition $\lambda$.  It is defined by the formula 
\[ w_\lambda(i) \, = \, 
\begin{cases} i + \lambda_{1-i} \; & \text{if} \;\; i \leq 0; \\
i - \lambda'_i \; & \text{if} \;\; i > 0.
\end{cases}\]
These are exactly the permutations such that $w(i) < w(i+1)$ for all $i \neq 0$.  In 
fact, $w_\lambda = w(\tau)$, where $\tau$ has $p_i = 0$ for all $i$.  The sequence 
$l_\bull$ is constructed so the
corners of $\lambda$ are at $(k_i,l_i)$; therefore  $q_i = l_i - k_i$.  

The \emph{multivariate Schur polynomial} of $\lambda$ is
the Schubert polynomial of the vexillary permutation $w_\lambda$.  
By the general vexillary formula, it has a determinantal formula:
\begin{cor}\label{multivar}
\[ 
\enS_{w_\lambda} = S_\lambda(c(1), \ldots, c(k_s)) \, = \, 
\sum_{\mu \subset \lambda}S_{\lambda/\mu}(a(1), \ldots, a(k_s)) \, S_\mu(c) ,
\]
where $c(k) = c \cdot a(k)$ and $a(k) = a(0,q_i)$ for $i$ minimal such that  $k_i \geq k$.  
\end{cor}
In particular, this is a polynomial in $\LLambda[y]$, and when the $y$ variables are set equal to $0$, it becomes $S_\lambda(c)$.  For example, for the partitions $(1)$, $(2)$, and $(1,1)$, these polynomials are $c_1$,  $c_2 + y_1 c_1$, and $c_1^2 - c_2 - y_0 c_1$.  These polynomials, with some of their sources and properties, are discussed in \cite[\S4 and App.~A]{LLS1}.

\section{Cauchy-Interpolation, Decomposition Formula}

The Schubert polynomials $\enS_w$, as $w$ varies over $\Sgp_\Z$, form a basis 
for $\LLambda[x,y]$ over  $\Z[y]$.   Interpolation gives a formula for writing a 
polynomial in terms of this basis.
Let $\eta \colon \LLambda[x,y] \to \Z[y]$ be the $\Z[y]$-algebra homomorphism obtained 
by mapping each $c_i$ to $0$ and each $x_i$ to $-y_i$.

\begin{prop}  
Any element of $\LLambda[x,y]$, has a unique expression as $f = \sum_{w \in \Sgp_\Z} a_w(y) \enS_w$, with  $a_w(y) = \eta(\del_w^x(f))$.
\end{prop}

For example, $x_1 = \enS_{s_1} - \enS_{s_0} - y_1\enS_{\id}$.  In the proposition, as usual, 
the notation $\del_w^x$ denotes $\del_{i_1}^x \circ \cdots \circ \del_{i_l}^x$ for any sequence of integers such that $w = s_{i_1}  \cdots  s_{i_l}$ and $l = \ell(w)$.

There is a useful \emph{decomposition formula} for these Schubert polynomials:

\begin{prop}\label{decomp}  Let $a$ and $b$ be Chern series,  
and let  $x$, $y$, and $t$ be sequences of variables $x_i$, $y_i$, and $t_i$  for 
all integers $i$.  Then, for any $w$ in $\Sgp_\Z$, 
\[ 
\enS_w(a\cdot b,x,y) = \sum_{v \cdot u \stackrel{.}{=} w} \enS_u(a,x,t)\, \enS_v(b,-t,y).
\]
\end{prop}

\begin{cor}\label{decompcor} For Chern series  $a$, $b$ and $c$,  sequences  $x$, $y$, $s$ and $t$ of 
variables, and $w$ in $\Sgp_\Z$,
\[
\enS_w(a\cdot c \cdot b,x,y) =  \sum_{v \cdot t \cdot u \stackrel{.}{=} w} \enS_u(a,x,s)\, \enS_t(c,-s,t) \,\enS_v(b,-t,y) .
\]
\end{cor}

The decomposition formula can be used in place of some Hopf algebra arguments and coproduct formulas in \cite{LLS1}.

\section{Stanley Polynomials, Formula for Schubert Polynomials}

The \emph{Stanley polynomial} $F_w$  of  $w$ in $\Sgp_\Z$  is the polynomial in 
$\LLambda = Z[c]$  defined by setting 
all the  $x$ and $y$ variables equal to $0$:
\[  F_w = \enS_w(c,0,0). \]

For examples,  $F_{s_k} = c_1$ for all $k$; $F_{3\,1\,2} = c_2$; $F_{2\,3\,1} = S_{(1,1)}(c) =
c_1^2-c_2$;  $F_{3\,2\,1} = S_{(2,1)}(c) = c_2 c_1 - c_3$; and $F_{2\,1\,4\,3} = c_1^2$.

\begin{prop} For any vexillary $w = w(\tau)$, 
 \[ F_{w} = S_{\lambda(\tau)}  .\]
\end{prop}

\begin{prop}  For any $w$ in $\Sgp_\Z$ and integer  $m$, 
\[  F_{\gamma^m(w)} = F_w . \]
\end{prop}

Stanley polynomials are known to be positive combinations of 
Schur polynomials:
\[
F_w = \sum j_\lambda^w \, S_\lambda(c).
\]
Fomin, Greene, Reiner and Shimozono \cite{FGRS} showed that 
 $j_\lambda^w$ is the number of tableaux of shape 
$\lambda$, strictly increasing in rows and columns, whose 
row reading word (read from top to bottom, in rows from 
right to left), is a reduced word  for 
$w$.  A bumpless pipedream formula is given in \cite{LLS1}.

The involution $\omega$ defined at the end of the first section restricts to an isomorphism between  $\Sgp_{-}$ 
and $\Sgp_{+}$, and determines an isomorphism $\omega \colon \Z[x_{+},y_{+}] \to \Z[x_{-},y_{-}]$, where the second ring is the polynomial ring in the variables $x_i$ and $y_j$  for  $i \leq 0$ and $j \leq 0$.  

For $u$ in $\Sgp_{-}$, $\omega(u)$  in $\Sgp_{+}$  has a double Schubert polynomial 
$\frak{S}_{\omega(u)}(x,y)$ in $\Z[x_{+},y_{+}]$.  The \emph{double Schubert polynomial} 
$\frak{S}_u(x,y)$ is defined to be the image of this by $\omega$, in  $\Z[x_{-}, y_{-}]$:
\[  \frak{S}_u(x,y) = \omega( \frak{S}_{\omega(u)}(x,y) ) .\]

For $w$ in $\Sgp_{\neq 0}$, write $w = u \cdot v$, with $u \in \Sgp_{-}$ and $v \in \Sgp_{+}$, and define 
the \emph{double Schubert polynomial} $\frak{S}_w(x,y)$ in $\Z[x,y]$ to be their product:
\[
  \frak{S}_w(x,y) = \frak{S}_u(x,y) \cdot \frak{S}_v(x,y) = \omega( \frak{S}_{\omega(u)}(x,y)) \cdot \frak{S}_v(x,y).
\]
(If the $y$ variables are set equal to $0$, these become the single polynomials $\frak{S}_w(x)$ defined in \cite{LLS1}.)
\begin{prop}  For any $w$ in $\Sgp_\Z$,  
\[  \enS_w(0,x,y) = \begin{cases} 
\frak{S}_w(x,-y)  & \text{ if }  w  \in \Sgp_{\neq 0} 
\\
0 & \text{ if }   w  \notin  \Sgp_{\neq 0} .
\end{cases} \]
\end{prop}

From Corollary~\ref{decompcor}, when the variables $a_i$, $b_i$, $s_i$ and $t_i$ are specialized to $0$, one recovers \cite[Cor.~4.5]{LLS1}:

\begin{prop}  
For any $w$ in $\Sgp_\Z$, 
\[
\enS_w = \sum_{\substack {v^{-1} \cdot t \cdot u \stackrel{.}{=} w \\
u, v \in \Sgp_{\neq 0}} }
\frak{S}_u(x)\, \frak{S}_v(y)\, F_t  .
\]
\end{prop}

\section{Products, Chevalley-Monk, and Transition}

For any $u$, $v$ and $w$ in $\Sgp_\Z$, there is a unique $c_{u\,v}^w$ in $\Z[y]$, homogeneous 
of degree $\ell(u) + \ell(v) - \ell(w)$, such that 
\[  \enS_u \cdot \enS_v = \sum c_{u\,v}^w \, \enS_w.\]

\begin{prop}  
\begin{enumerate}[wide,itemsep=5pt]
\item Each $c_{u\,v}^w$ is in $\Z_{\geq 0}[ \dots, y_i - y_{i+1}, \ldots]$, i.e., it is 
a nonnegative linear combination of products of the $y_i - y_{i+1}$, $i \in \Z$.  

\item $c_{u\,v}^w = 0$  unless $u$ and $v$ are contained in $w$ in the Bruhat order.  

\item  For all $u$, $v$, $w$ in $\Sgp_\Z$, and for all integers $m$,
\[  c_{\gamma^m(u) \, \gamma^m(v)}^{\gamma^m(w)} = \gamma^m(c_{u\,v}^w); \]
note that $\gamma^m(y_i - y_{i+1}) = y_{m+i}- y_{m+i+1}$.  

\item For $u$, $v$ and $w$ in $\Sgp_{+}$, $c_{u\,v}^w$  is the same 
structure constant as for double Schubert polynomials:
\[ \frak{S}_u(x,-y) \frak{S}_v(x,-y) \, = \, \sum  c_{u\,v}^w \frak{S}_w(x,-y) .\]
\end{enumerate}
\end{prop}

\noindent
(For single Schubert polynomials, i.e., $y=0$, see \cite[Prop.~3.20]{LLS1}.)

The \emph{Chevalley-Monk} formula also extends to this setting.  For any 
integers $i < j$, let $t_{i\,j}$ be the transposition that interchanges $i$ and $j$. 
Recall that $w \, t_{i\,j}$ is a \emph{cover} of $w$ in the Bruhat order, i.e.
$\ell(w t_{i\,j}) = \ell(w) + 1$, exactly when $w(i) < w(j)$ and there is no $i < l < j$ 
with $w(i) < w(l) < w(j)$; this covering relation is denoted $w \, \nearr \, w\, t_{i\,j}$.

\begin{prop}  For all $w \in \Sgp_\Z$ and integers $k$,
\[  \enS_{s_k} \enS_w = \sum_{\substack {i \leq k < j \\ w \nearr w \, t_{i\,j}}}   
\enS_{w \,t_{i,j}} + \sum_{i \leq k} (y_i - y_{w(i)}) \enS_w .\]
\end{prop}

\begin{cor} For  $w \in \Sgp_\Z$ and  $k \in \Z$,
\[ (x_k + y_{w(k)}) \enS_w = \sum_{\substack {k < j \\ w \nearr w \, t_{k\,j}}}   
\enS_{w \,t_{k,j}} - \sum_{\substack {i < k \\ w \nearr w \, t_{i\,k}}} \enS_{w \,t_{i,k}}   . \]
\end{cor}

 Note that even if $u$ and $v$ are in $\Sgp_{+}$, there can 
be $w$ not in $\Sgp_{+}$ for which $\enS_w$ appears in the product of $\enS_u$ and $\enS_v$.  
For example,  $\enS_{2\,1} \enS_{2\,1} = \enS_{3\,1\,2} + \enS_{1\, 2 \,0} + (y_1 - y_2) \enS_{2\,1}$.

The following \emph{transition formula} is often an efficient way to 
calculate Schubert polynomials, or for inductive proofs.

\begin{prop}
Let $w$ be in $\Sgp_\Z$, and $r < s$ a pair of integers  such that $w(r) > w(s)$ and no
$r < q < s$ has $w(q)>w(s)$ and no $q > s$ has $w(r) > w(q) > w(s)$.  Let $v = w\,t_{r\,s}$. 
Then 
\[  
\enS_w \,=\, (x_r + y_{w(s)}) \enS_v + \sum_{\substack {i < r \\ v \nearr v \, t_{i\,r}}} \enS_{v \,t_{i,r}}  . 
\]
\end{prop}

Unlike in the classical setting, cf.~\cite{Bi1}, applying this repeatedly does not terminate, as there may always 
be descents.  However, transition gives an algorithm to write $\enS_w$ as a sum of products 
of linear factors $x_i + y_j$  times Schubert polynomials $\enS_v$, with $v$ vexillary, for which we have explicit determinantal formulas.  In fact, one needs only $v$ for which some $\gamma^m(v)$ is a dominant permutation in $\Sgp_{+}$.

For example, with $\alpha_{i\,j} = x_i + y_j$,
\begin{align*}
  \enS_{1\,0\,3\,2} &= (x_2+y_2) \enS_{1\,0\,2\,3} + \enS_{1\,2\,0\,3} + \enS_{2\,0\,1\,3} \\
  &=
\alpha_{2\,2} c_1+(c_1^2-c_2+x_1c_1)+(c_2+y_1c_1) = c_1^2+(\alpha_{1\,1}+\alpha_{2\,2})c_1.
\end{align*}

\section{Localization}

For any $v$ in $\Sgp_\Z$, there is a \emph{localization homomorphism} 
\[
 \phi_v \colon \LLambda[x,y] \to \Z[y] , 
\] 
a homomorphism of $\Z[y]$-algebras, defined by sending $x_i$ to $-y_{v(i)}$, and sending $c$ to 
\[
  \prod_{i \in  \Z_{\leq 0} \, \cap \, w(\Z_{> 0})} (1+y_i) \,\, / \prod_{j \in  \Z_{>0} \, \cap  \,w(\Z_{\leq 0})} (1+y_j) .
\]

Let $F(\Sgp_\Z,\Q(y))$ be the $\Q(y)$-algebra of all functions from the set $\Sgp_\Z$ to the 
quotient field $\Q(y)$ of $\Z[y]$.  There is a homomorphism 
\[
 \Phi \colon \LLambda[x,y] \to F(\Sgp_\Z,\Q(y)), \;\;\; P \mapsto [v \mapsto \phi_v(P)],
\]
which is an embedding. For any integer $i$,  define the endomorphism $A_i$ of $F(\Sgp_\Z,\Q(y))$  by the formula 
$ A_i(f)(v) = \frac{f(s_i \, v) - f(v)}{y_{v(i)} - y_{v(i+1)}}$.  Then, as in \cite{Ar},
$A_i \circ  \Phi = \Phi \circ \del_i^x $ for all integers $i$.

There is also a straight-forward generalization of the AJS-Billey formula \cite{AJS}, \cite{Bi2}, which gives an explicit expression for  $\phi_v(\enS_w)$, for any $v$ and $w$ in $\Sgp_\Z$, as a sum of products of $y_i - y_j$ for $i < j$. 

\section{Twisted Schubert Polynomials} \label{twisted} 

It is useful to have \textit{twists} of Schubert polynomials.  
These correspond to the geometric notion of tensoring (twisting) by a line bundle.  
There is a new variable $z$, corresponding to the first 
Chern class of the line bundle, which becomes part of the algebra. 
In types C and D, this twist is crucial, allowing one to have 
symplectic and quadratic forms with values in a nontrivial line bundle.

In type A, twists are algebraically simpler.  The algebra $\LLambda$ is replaced by $\bLambda = \LLambda[z] = \Z[z,c]$, the polynomial ring in variables $z$ of degree $1$ and $c_k$ of degree $k$.  We define \emph{twisted enriched Schubert polynomials} 
$\tenS_w = \tenS_w(c,x,y,z)$ in $\bLambda[x,y]$, for any $w$ in $\Sgp_\Z$.  The simplest way to do this is to replace each variable $x_i$ that appears by $x_i - z$:
\[
 \tenS_w = \tenS_w(c,x,y) = \tenS_w(c,x,y,z) := \enS_w(c,x-z,y) .
\]
For example, setting $\balpha_{i j} = x_i + y_j - z$, 
\[
 \tenS_{s_k} = \begin{cases} c_1 + \sum_{i=1}^k \balpha_{i\,i} \text{ for } k \geq 0 \\
c_1 - \sum_{i=k+1}^0 \balpha_{i\,i} \text{ for } k < 0.
\end{cases} 
\]
%

There is also a formula for these polynomials that subtracts $z$ from the $y$ variables, provided the $c$ variables are appropriately modified:

\begin{prop}
\[
 \tenS_w(c,x,y) = \enS_w(\theta(c),x,y-z),
\]
where $\theta(c)_k = \sum_{i=1}^k \binom{k-1}{i-1} z^{k-i} c_i$. 
\end{prop}

Setting $z$ to $0$ in $\tenS_w$ recovers the polynomials $\enS_w$ we have been 
studying.   Most of the identities stated here extend automatically to this twisted setting, 
by replacing each occurrence of $x_i$ by $x_i - z$.  For example, there are difference 
operators $\del_i^x$ and $\del_i^y$ on these rings, using the same formulas, except in 
the formula for $\del_0^x(S_\lambda(c))$ of Proposition~\ref{del0}, one needs to replace the 
$x_1$ by $x_1 - z$ and the $-x_0$ by $z - x_0$.  In localization formulas, one replaces 
$\alpha_{i\,j} = x_i + y_j$ by $\balpha_{i\,j} = x_i + y_j - z$.

Another formula that needs modification is the duality formula.  The involution $\omega$ of $\LLambda$ lifts to an involution, denoted $\tomega$, that takes $c_k$ to
\[
\tomega(c)_k = \tomega(c_k) = 
\sum_{i=1}^k \tbinom{k-1}{i-1}(-z)^{k-i}S_{1^i}(c).
\]

\begin{prop}  \[
\tenS_{w^{-1}}(c,x,y) = \tenS_w(\tomega(c) ,y,x).
\]\label{tilde}
\end{prop}

The formula $\tenS_{\omega(w)} = \tomega(\tenS_w)$ requires defining $\tomega(x_i) = z - x_{1-i}$ and
$\tomega(y_j) = z - y_{1-j}$. 

The formulas for products are unchanged: the $c_{u\,v}^w$ are the same, and do not 
involve the $z$ variable.

The enriched Schubert polynomials also satisfy an \textit{invariance property}, generalizing the fact that $\frak{S}_w(x+v,y+v) = \frak{S}_w(x,y)$ for any variable $v$.  Here one has
 \[  \tenS_w(c,x+v,y,z+v) = \tenS_w(c,x,y,z).\]  
Equivalently, if $\theta_v(c)$ is defined by the formula 
$\theta_v(c)_k = \sum \binom{k-1}{i-1} v^{k-i} c_i$, then
\[ \tenS_w(\theta_v(c),x,y+v,z+v) = \tenS_w(c,x,y,z).\]

\section{Degeneracy loci}

We are given vector bundles $V$ and $W$, and a line bundle $L$, on a variety $X$, assumed to be nonsingular for simplicity.  We are given a bilinear mapping $V \times W \to L$, and we are given two flags of subbundles of $V$, labeled as follows:
\[
0 \subset E_{n} \subset E_{n-1} \subset \dots \subset E_{1} \subset E_{0}  \subset E_{-1} \subset \dots \subset E_{-m} \subset  V ,
\]
\[
 0 \subset F_{n} \subset F_{n-1} \subset \dots \subset F_{1} \subset F_{0}  \subset F_{-1} \subset \dots \subset F_{-m} \subset  W .
\]
The key assumption is that $E = E_0$ and $F= F_0$ have the same rank, 
so $\rank(E_i) = \rank(E) -i  =  \rank(F_i)  = \rank(F) - i$ for all $i$.

Let $w$ be in $\Sgp_\Z$, and assume  $m$ and $n$ are large enough so all integers moved by $w$ lie in the interval $(-m,n]$.  
Define 
\[k_w(p,q) = \#\{a \leq p \mid w(a) > q \}.\]
We 
have the \emph{degeneracy locus} $\Omega_w = \Omega_w(E_\bull,F_\bull)$ in $X$, 
defined to be the locus where the nullity of each  $E_p \otimes F_q \to L$ is at least  $k_w(p,q)$ for all $p$ and $q$.  That is,
\[
 \Omega_w(E_\bull,F_\bull) \; : \;\; \dim \Ker (F_q  \to \Hom(E_p,L)) \geq k_w(p,q) \;\; \forall \; p,q .
\]
This is a subscheme of $X$ locally defined by a collection of determinants.  An important particular case is when $w = w(\tau)$ is vexillary, coming from a triple $\tau = (k_\bull,p_\bull,q_\bull)$.  Then the locus $\Omega_w(E_\bull,F_\bull)$ is the locus $\Omega_\tau(E_\bull,F_\bull)$ where the dimension of the kernel of the map from $F_{q_i}$ to $\Hom(E_{p_i},L)$  is at least  $k_i$, for each $i$.

We set
\[
c = c(E^*\otimes L - F) ,\;\; x_i = c_1(E_{i-1}/E_i), \;\;  y_j = c_1(F_{j-1}/F_j).
\]

\begin{thm}\label{degloci}
For any $w$, the codimension of 
$\Omega_w(E_\bull,F_\bull)$ in $X$ is at most the length $\ell(w)$. If the  codimension  is equal to $\ell(w)$, its class $[\Omega_w(E_\bull,F_\bull)]$ is given by the polynomial $\tenS_w(c,x,y)$.
\end{thm}

\begin{cor}
For any type A triple $\tau = (k_\bull,p_\bull,q_\bull)$, the codimension of $\Omega_\tau(E_\bull,F_\bull)$ is at most $|\lambda|$, where $\lambda = \lambda(\tau)$, and, when equal, its class is given by the polynomial 
\[
 S_{\lambda}(\mathbf{c}(1), \ldots, \mathbf{c}(k_s)) = \sum_{\mu \subset \lambda} S_{\lambda/\mu}(\mathbf{a}(1), \ldots, \mathbf{a}(k_s)) S_\mu(c),
\]
where $\mathbf{a}(k)$ is defined as in \S4, but replacing each  $x_i$ in $a(k)$ by $x_i - z$.
\end{cor}

If the loci do not have the predicted dimensions, or the ambient variety is singular, there are  \textit{refined} (or \textit{virtual}) classes, as in Theorem 14.3 of \cite{F1}.

Consider the \emph{split} case, where $V = \bigoplus_{i = -m}^n L_i$ and  $W = \bigoplus_{i = -m}^n M_i$ are direct sums of line bundles  $L_i$ and $M_j$  with first 
Chern classes $x_i$ and $y_j$.  Set  $E_p = \bigoplus_{i > p} L_i$, 
$F_q = \bigoplus_{j > q} M_j$.  The bilinear form is given by a square matrix $A = (a_{i\,j})$, with rows and columns numbered from $-m$ to $n$, where 
$a_{i\,j}$ is a section of $\Hom(L_i \otimes M_j,L)$.  The locus $\Omega_w(E_\bull,F_\bull)$ requires each \textit{southeast} subrectangle of $A$ consisting of rows strictly below row $p$ and columns strictly right of column $q$ to have nullity at least $k_w(p,q)$.  Its class is given by $\tenS_w(c,x,y)$, where now $c$ specializes to $\prod_{i = 1}^n  \frac{(1 + (z - x_i)t)}{(1+y_i t)}$. 

\begin{prop}
The polynomial $\tenS_w$ is uniquely determined by these formulas for degeneracy loci in the split case.
\end{prop} 

If one takes a square matrix $A = (a_{i\,j})$ with indeterminate coefficients, and  gives $a_{i\,j}$ the \emph{weight} $z - x_i - y_j$, the \emph{matrix determinantal loci} $\Omega_w$ determine these Schubert polynomials, much as the loci studied by Knutson and Miller \cite{KM} did for classical Schubert polynomials.

\medskip

There is an alternative (equivalent) degeneracy locus setup, which is often useful.  Here one has a vector bundle $W$, with two flags of subbundles $G_\bull$ and $F_\bull$ of $W$
\[
 \cdots \subset G_{-1} \subset G_0 = G \subset G_1 \subset G_2 \subset \cdots \subset W,
\]
\[
 \cdots \subset F_2 \subset F_1 \subset F_0 = F \subset F_{-1} \subset \cdots \subset W,
\]
indexed so that $\rank G_p = \rank(G)+p$ and $\rank F_q = \rank(F)-q$, and now
\[ \rank(G) + \rank(F) \, = \, \rank(W).\]   Now we have the \text{intersection locus} $ \Omega_w(G_\bull \cap F_\bull)$ defined by the conditions
\[
 \Omega_w(G_\bull \cap F_\bull)\;  : \;\; \dim( G_p \cap F_q ) \geq k_w(p,q) \;\; \forall \; p,q .
\]

\begin{prop}\label{altdegloci}
The codimension of $\Omega_w(G_\bull \cap F_\bull)$ in $X$ is at most the length $\ell(w)$.  If the codimension  is equal to $\ell(w)$, its class $[\Omega_w(G_\bull \cap F_\bull)]$ is given by the polynomial $\enS_w(c,x,y)$, where $c$ maps to $c(W-G-F)$, $x_i$ maps to $- c_1(G_i/G_{i-1})$, and $y_j$ maps to $ c_1(F_{j-1}/F_j)$.
\end{prop}

One can deduce formulas due to Knutson and Pawlowski for {\em graph Schubert varieties} \cite{P}.  On a variety $X$, we have vector bundles $U$ and $V$, both of rank $n$, with flags of subbundles
\[
G_\bull : G_1 \subset G_2 \subset \cdots \subset V \quad \text{and} \quad F_\bull: U = F_0 \supset F_1 \supset F_2 \supset \cdots .
\]
Let $W=U\oplus V$, and form the Grassmann bundle $\bGr(n,W) \to X$.  Let $S \subset V$ be the tautological rank $n$ subbundle.  For a permutation $w\in \Sgp_n$, the \define{graph locus} $\Gamma\Omega_w \subset \bGr(n,V)$ is defined by the conditions
\[
 \dim( S \cap F_q \oplus G_p ) \geq k_w(p,q) \quad \forall \; p,q >0.
\]

\begin{cor}[cf.~{\cite[Theorem~3.11]{P}}]\label{Pawlowski}
If $\Gamma\Omega_w$ has codimension $\ell(w)$, then its class is
\[
  [\Gamma\Omega_w] = \enS_w(c,x,y),
\]
where $c=c(U-S)$, $x_i=-c_1(G_i/G_{i-1})$, and $y_j=c_1(F_{j-1}/F_j)$.
\end{cor}

Pawlowski considers the case where $X = \Fl(\C^n) \times \Fl(\C^n)$, the bundles $U$ and $V$ are trivial, and $G_\bull$ and $F_\bull$ are tautological flags from the two factors.  In this setting, he shows that the locus $\Gamma\Omega_w \subseteq \Gr(n,\C^{2n})\times \Fl(\C^n)\times\Fl(\C^n)$ is irreducible of codimension $\ell(w)$, with class $\enS_w(c,x,y)$.

\section{Type C}

We sketch the type C story, to illustrate the close analogy with the type A story, and to set up what is needed to prove our theorem comparing Schubert polynomials in the two types.  Details can be found in \cite{AF1} and \cite{AF4}.

The basic ring $\bGamma$  for type C is a quotient ring of $\bLambda = \Z[z,c] = \Z[z,c_1,c_2, \ldots ]$ modulo 
an ideal of relations.  Define, for any nonnegative integers $p \geq q$, 
\[ 
 C_{p\,q} = \sum_{0 \leq i \leq j \leq q} (-1)^j \, (\tbinom{j}{i}+\tbinom{j-1}{i}) \, z^i \, c_{p+j-i} \, c_{q-j} .
\]
Set
\[ 
\bGamma = \Z[z,c]/(C_{1\,1},C_{2\,2}, C_{3\,3}, \ldots ) .
\]

These relations can be motivated by geometry.  Suppose $E$ and $F$ are two vector bundles of the same rank, and $L$ is a line bundle; set $c_k = c_k(E^*\otimes L - F)$, and $z = c_1(L)$, then the relations say that $c_k(F^*\otimes L - E) = c_k$ for all $k$, i.e., that $\tomega(c) = c$.  In the type C setting, $E$ and $F$ will be maximal isotropic subbundles of a vector bundle $V$ with an $L$-valued symplectic form, which implies that  $\Hom(E,L) = V/E$ and  $\Hom(F,L) = V/F$, and
\[ 
c(E^*\otimes L - F) = c(V - E - F) = c(F^*\otimes L - E) .
\]

Note that sending $c_i$ to $q_i$ determines an isomorphism $\bGamma/ z \, \bGamma \cong \GGamma$, where $\GGamma$ is the classical ring generated by the Schur Q-functions $q_1, q_2, \ldots$.  

This $\Z[z]$-algebra has a basis of \textit{pfaffians} $Q_\lambda(c)$, as $\lambda$ 
varies over all strict partitions.  For $\lambda = (p)$, $Q_\lambda(c) $ is $c_p$; for $\lambda = (p,q)$,   $Q_\lambda(c)$ is $C_{p\, q}$.  In general, $Q_\lambda(c)$ is the pfaffian of the matrix whose $(i,j)$ entry is $C_{\lambda_i \, \lambda_j}$, adding a $0$ to the end of $\lambda$ if its length is odd, so that the matrix has an even number of rows and columns.

The \emph{type C twisted double Schubert polynomial} $\tenS_w^C = \tenS_w^C(c,x,y)$ is defined for any signed permutation $w$; 
it lives in $\bGamma[x_+,y_+]$.  The construction can be given as follows.

For any strict partition $\lambda$ of length $r$, and any Chern series $c(1), \ldots, c(r)$, denote by $ \Pf_\lambda(c(1), \ldots, c(r))$
the pfaffian of the alternating matrix whose $(k,l)$ entry, for $1 \leq k < l \leq r$, is 
\[ 
\sum_{0 \leq i \leq j \leq \lambda_l} (-1)^j \,(\tbinom{j}{i}+\tbinom{j-1}{i}) \, z^i \, c(k)_{\lambda_k - \lambda_l +j-i}\,  c(l)_{\lambda_l-j} .
\]
(As before, append a $0$ to $\lambda$ if its length is odd.)

We start with those $w$ we call vexillary, whose Schubert polynomials are given by explicit pfaffians \cite{af-vex} (see also \cite{BL}).  These come from a \emph{type C triple} $\tau = (k_\bull,p_\bull,q_\bull)$ of some 
length $s$.  Here $k_\bull$ is a sequence of $s$ strictly increasing positive integers, and $p_\bull$ and $q_\bull$ are both 
sequences of $s$ weakly decreasing positive integers: 
\[ 
 k_1 < k_2 < \dots < k_s, \;\;\; p_1 \geq p_2 \geq \dots \geq p_s,  \;\;\; q_1 \geq q_2 \geq \dots \geq q_s. 
\]
These are required to satisfy the condition
\[ 
(p_i - p_{i+1}) + (q_i - q_{i+1}) \, > \, k_{i+1} - k_i \tag{*} 
\] 
for $1 \leq i < s$.  

A triple $\tau$ determines a strict partition $\lambda = \lambda(\tau)$ of length $k_s$; it is the minimal strict partition with $\lambda_{k_i} = p_i + q_i - 1$; that is,  for $1 \leq k \leq k_s$, take $i$ minimal so $k_i \geq k$, and set $\lambda_k = p_i+q_i - 1+k_i - k$.

A triple $\tau$ also determines a signed permutation $w = w(\tau)$, which is the signed permutation of minimal length 
such that 
\[
  \# \{ a \geq p_i \mid -w(a) \geq q_i \} \, = \, k_i 
\]
for $1 \leq i \leq s$.  

Set 
\[ 
 c(k) = c \cdot a(k),  \text{ with } a(k) = \prod_{a < p_i} (1 + x_a) \cdot \prod_{b < q_i} (1+y_b) ,
\]
where again $i$ is minimal so $k_i \geq k$.  The Schubert polynomial of 
$w(\tau)$ is given by the formula
\[
\begin{split}  \tenS_{w(\tau)}^C =& \Pf_{\lambda(\tau)}(c(1), \ldots, c(k_s)) \\
=& \sum \Pf_{\lambda/\mu}(a(1), \ldots, a(k_s)) Q_\mu(c),
\end{split}
\]
the sum over strict partitions $\mu$ contained in $\lambda$.  A positive formula for these coefficients can be derived using the method of non-intersecting paths.
There are difference operators $ \del_i^x$  and $\del_i^y$ acting on the rings $\bGamma[x,y]$, which send the twisted Schubert polynomials to each other, so that all can be constructed from these pfaffians.  In fact, these Schubert polynomials satisfy analogues of all the properties described in the previous sections for the type A Schubert polynomials; details can be found in \cite{AF1} and \cite{AF4}.

We will need the corresponding type C degeneracy loci formula.  We start with a vector bundle $V$ of even rank on a (nonsingular) variety $X$, with an \emph{ L-valued symplectic 
form}.  This is a non-degenerate, alternating, bilinear form from 
$V \times V$ to a line bundle $L$.  Assume  maximal isotropic  subbundles
$E$ and $F$ of $V$ are specified 
(each of rank half the rank of $V$, with the symplectic form vanishing on each of them).   
We are given flags $E = E_1 \supset E_2 \supset \dots $, 
and $F = F_1 \supset F_2 \supset \dots $, where now 
$\rank(E_i) = \rank(F_i) = \rank(E)+1-i$.  These are completed to complete flags 
\[
0 \subset   \dots \subset E_2 \subset E_1 = E \subset E_{-1} \subset E_{-2} \subset \dots  \subset V ,
\]
by setting $E_i = E_{1-i}^\perp$ for $i < 0$, so $\rank(E_i) = \rank(E)-i$ for $i < 0$.  Similarly  $F_i = F_{1-i}^\perp$, 
giving a complete flag for them as well.  For any signed permutation $w$ of $\{\pm 1, 2, \ldots, \pm \rank(E)\}$,
we have the degeneracy locus  $\Omega_w^C(E_\bull,F_\bull)$
defined by setting 
\[  
\Omega_w^C(E_\bull,F_\bull)  \; : \;\; \dim(E_p \cap F_q) \geq \#\{ a \geq p \mid -w(a) \geq q \} ,
\]
for $p \geq 1$ and arbitrary $q$.  
For a type C triple $\tau$ of length $s$, we have the locus 
\[ 
 \Omega_\tau^C(E_\bull,F_\bull)  \; : \;\; \dim(E_{p_i} \cap F_{q_i}) \geq k_i 
\]
for $1 \leq i \leq s$.  Note that  $\Omega_{w(\tau)}^C(E_\bull,F_\bull) =  \Omega_\tau^C(E_\bull,F_\bull)$.

The type C Schubert polynomials give, and are determined by, formulas for these degeneracy loci.  For this, set 
$c = c(V - E - F)$, $z = c_1(L)$,  $x_i = c_1(E_{i-1}/E_i)$ and  $y_i = c_1(F_{i-1}/F_i)$.  Note that  $c_1(E_{-i}/E_{1-i}) = z - x_i$ and $c_1(F_{-i}/F_{1-i}) = z - y_i$ for $i \geq 1$.

\begin{prop}
The codimension of $\Omega_w^C(E_\bull,F_\bull)$ in $X$ is at most $\ell(w)$.  When equality holds, 
\[
 [ \Omega_w^C(E_\bull,F_\bull) ] = \tenS_w^C(c,x,y) . 
\]
\end{prop}

The fact that 
$ \Omega_{w^{-1}}^C(E_\bull,F_\bull) = \Omega_w^C(F_\bull,E_\bull) $
corresponds to the identity
\[ 
\tenS_{w^{-1}}^C(c,x,y) =  \tenS_w^C(c,y,x).
\]

The degeneracy locus formula implies a succinct formula for the locus $\Omega_k(\sigma)$ where a symmetric map of vector bundles $\sigma \colon E \to E^* \otimes L$ drops rank by at least $k$, for $k \leq \rank(E)$, for a  vector bundle $E$ on a (nonsingular) variety $X$:
\[  \Omega_k(\sigma) \, : \; \dim (\Ker ( E \to  E^* \otimes L )) \geq k .\]
For this, set $V = E \oplus \Hom(E,L)$.  This has a standard $L$-valued symplectic form: $\langle (v, \phi), (v',\phi') \rangle = \phi(v') - \phi'(v)$.  The fact that $\sigma$ is symmetric implies that the graph of $\sigma$ is a maximal isotropic subbundle of $V$.   Apply the degeneracy locus formula to this graph and the graph $E \oplus 0$ of the zero morphism, noting that $\Omega_k(\sigma)$ is the locus where these two maximal isotropic subbundles meet in dimension at least $k$:   

\begin{cor} 
The codimension of  $\Omega_k(\sigma)$ in $X$ is at most $\binom{k+1}{2}$.  If equality holds, then its class is
given by a pfaffian in the variables $c_k$:
\[ 
 [\Omega_r(\sigma) ]  \, = \, Q_{\lambda}(c), 
\]
where $\lambda = (k, k-1, \ldots, 2,1)$, and $c$ maps to $c(E^*\otimes L - E)$.
\end{cor}

\noindent
In the case where $L$ is trivial, this is equivalent to a formula proved by Pragacz \cite[Proposition~7.8]{pragacz88}.

\section{Type A to Type C}

In this section we write, for $w$ in $\Sgp_{+}$, $\tenS_w^A = \tenS_w^A(c,x,y)$ 
for the type A twisted enriched polynomial discussed in this paper; it is a polynomial in
$\bLambda[x_+,y_+]$, the polynomial ring in the positive $x$ and $y$ variables.

There is a canonical homomorphism of $\Z[z]$-algebras 
\[ \bLambda[x_+,y_+] \to \bGamma[x_+,y_+] \]
sending $c_k$ to $c_k$, $x_i$ to $x_i$, $y_i$ to $y_i$.  In fact, 
if $\tomega(c)$ is defined as in Section~\ref{twisted}, then 
$c$ and $\tomega(c)$ have the same image.\footnote{This specializes to the 
classical fact that the map from $\LLambda$ to $\GGamma$ identifies the 
symmetric functions $h_i$ and $e_i$.} From this 
point of view the following theorem is natural, although it is not 
obvious from the algebraic constructions of the polynomials.

\begin{thm}\label{AtoC}  For $w$ in $\Sgp_{+}$, this canonical homomorphism 
sends  $\tenS_w^A$ to  $\tenS_w^C$.
\end{thm}

\begin{ex}  For $w = {3\,2\,1}$, it follows from the identities 
in $\bGamma$ that
$c_1^2 = 2 c_2 + z c_1$ and $c_2 c_1=c_{2\,1} + 2c_3 + z c_2$.  Therefore
$\tenS_w^A$ maps to $c_3 + c_{2\,1} + (x_1+y_1+\balpha_{1\,2}+\balpha_{2\,1})c_2
+ (y_1+\balpha_{1\,2})(x_1+\balpha_{2\,1})c_1 + \balpha_{1\,1}\balpha_{1\,2}+\balpha_{2\,1}$,
which is $\tenS_w^C$.
\end{ex}

\section{Intersection Rings}

The rings $\bLambda$ and $\bGamma$ control the intersection rings of Grassmann bundles in types A and C.  Write $A^*(X)$ for the Chow ring of a nonsingular variety $X$, but other intersection rings, such as singular cohomology in the complex case, can be used as well.   

In type A, suppose one has a non-degenerate bilinear mapping $V \times W \to L$ on a non-singular variety $Y$, and  a subbundle $F$ of $W$  on $Y$.  Let  $X \to Y$ be the Grassmann bundle of subspaces of fibers of $V$ of the same dimension as the rank of $F$.  (As usual, we use the same notation $V$, $W$, $L$, and $F$ for the pullbacks of these bundles to $X$.)  On $X$ we have also the tautological subbundle $E$ of $V$.  There is a canonical surjection 
\[  \bLambda \otimes_{\Z[z]} A^*(Y) \twoheadrightarrow A^*(X), \]
taking $c_k$ to $c_k(\Hom(E,L) - F)$.  In fact, the images of the Schur polynomials $S_\lambda(c)$ 
as $\lambda$ varies over partitions with length at most the rank of $F$ and entries at most the rank of $V/F$, give a basis for $A^*(X)$ over $A^*(Y)$.

In type C one has a symplectic form $V \times V \to L$ on $Y$, with a maximal isotropic subbundle $F$ of $V$.  Now $X \to Y$ is the bundle of maximal isotropic subspaces of fibers of $V$, with $E$ the tautological subbundle.  There is a canonical surjection 
\[  \bGamma \otimes_{\Z[z]} A^*(Y) \twoheadrightarrow A^*(X), \]
taking $c_k$ to $c_k(V - E - F)$.  The images of the $Q_\lambda(c)$, as $\lambda$ ranges over all strict partitions with entries at most the rank of $F$, give a basis for $A^*(X)$ over $A^*(Y)$.

There are similar presentations of the corresponding flag bundles, using the rings $\bLambda[x,y]$ and $\bGamma[x,y]$, and tensoring over $\Z[z,y]$.

\section{On Proofs}

The automorphisms $\theta$, $\omega$, and $\tomega$ used in this paper are 
general algebraic constructions, motivated by basic identities of Chern classes in geometry.  
For any commutative ring $R$, the polynomial ring $R[c,v] = R[c_1, c_2, \ldots][v]$  has an automorphism $\theta_v$ that takes $c_k$  to $\sum \binom{k-1}{i-1} v^{k-i} c_i$.  If $E$ and $F$ are vector bundles of the same rank on some variety, and $M$ is a line bundle, sending $v$ to the first Chern class of $M$, and sending $c$ to $c(E - F) = c(E)/c(F)$, then 
$\theta_v(c)$ maps to $c(E \otimes M^* - F \otimes M^*)$   Equivalently, 
\[ 
\theta_v(c(\xi)) =  c(\xi \cdot [M^*]) 
\]
for any $\xi$ of rank $0$ in the Grothendieck ring of vector bundles on the variety.
For another variable $u$,  $\theta_u$, $\theta_v$, and $\theta_{u+v}$ act on $R[c,u,v]$, 
and
\[ 
\theta_{u+v} = \theta_u \circ \theta_v ; 
\]
in particular  $\theta_{-v} = (\theta_v)^{-1}$.

Let $\omega$ be the involution of R[c]  that takes  $c$  to 
$1/(1 - c_1 + ...)$; so $\omega$ takes  $S_\lambda(c)$  to  $S_{\lambda'}(c)$, and $c(E - F)$ to $c(F^* - E^*)$, and $c(\xi)$ to $c(-\xi^*)$.  Then
\[ 
\omega \circ \theta_v = \theta_{-v} \circ \omega. 
\]
The lifted $\tomega$ on $\bLambda = \Lambda[z]$ is defined to be $\omega \circ \theta_z$.
So $\tomega \circ \theta_v = \theta_{-v} \circ \tomega$ for any  $v$.  And $\tomega(c(\xi)) = c(-\xi^* \cdot [L])$. 

Proposition~\ref{gamma} is proved in \cite[\S4]{AF2}. 

Theorem~\ref{DefThm} has can be deduced from several assertions here, e.g. from the degeneracy formula, or from the construction of these polynomials from vexillary polynomials and difference operators.  Or it can be proved directly, by verifying an appropriate supersymmetry.   It can also be deduced from \cite{LLS1}.  The fact that one need only check for a small number of $m$ follows from the following elementary fact:

\begin{lem}
For positive integers $n$, $p$, $q$, and indeterminates $c_1, \ldots, c_n$, $h_1, \ldots, h_p$, and $e_1, \ldots e_q$, 
the map $\Z[c_1,\ldots,c_n] \to \Z[h_1,\ldots, h_p, e_1, \ldots, e_q]$ that sends 
$c_k$ to $\sum_{i+j=k}h_i\,e_j$ is one-to-one if $n \leq p+q$.  
\end{lem}
\noindent (One can take $n = p+q$, $h_i$ to be the $i^{\text{th}}$ elementary symmetric polynomial in $p$ variables, and $e_j$ is the $j^{\text{th}}$ elementary symmetric polynomial in another set of $q$ variables; then $c_k$ is the $k^{\text{th}}$ elementary symmetric polynomial in both sets of variables, in which case it is obvious.)

Proposition~\ref{omega} has several proofs, besides being equivalent to Proposition 4.17 of \cite{LLS1}.  It also follows from the fact that 
\[ 
\del_i \circ \omega = \omega \circ \del_i 
\] 
for the $x$ and $y$ difference operators, and the characterization of the Schubert polynomials by Theorem~\ref{DiffOp}.  We will also see a proof via degeneracy loci.

Proposition~\ref{LLS1} follows from the identification of $c$ with the power series $\prod_{i \leq 0} \frac{1+y_i}{1-x_i}$, which becomes $\prod_{i \leq 0} \frac{1-a_i}{1-x_i}$ in the ring $\overleftarrow{R}(x;a) $.

Proposition~\ref{del0} is proved in \cite[\S4]{AF2}.  

Theorem~\ref{DiffOp} follows from the fact that no nonzero polynomial in $\LLambda[x,y]$ is annihilated by all the difference operators.  This is easily deduced from the fact that any element of $\LLambda$  annihilated by $\del_0$ must be an integer, which follows e.g.~from the explicit formula for $\del_0$ in Proposition \ref{del0}.

Corollary~\ref{vexcor} follows from the general identity in \cite[\S1]{AF2}.

The irreducibility of the vexillary Schubert polynomials in Corollary~\ref{irred} follows from the elementary and known fact that their images $S_\lambda(c)$ are irreducible in $\LLambda$, when the $x$ and $y$ variables are set equal to $0$.

Proposition~\ref{decomp} follows from the following identity of double Schubert polynomials, together with Theorem~\ref{DefThm}:
\[ 
\frak{S}_w(x,y) = \sum_{v \cdot u \stackrel{.}{=} w} \frak{S}_u(x,t) \, \frak{S}_v(t,y). 
\]

The formulas in the section on degeneracy loci uniquely determine the (twisted) enriched Schubert polynomials, so they can be used to prove many of the propositions in this paper. 
For example, the simple identity  
$\Omega_{w^{-1}}(E_\bull,F_\bull) = \Omega_w(F_\bull,E_\bull)$, together with the fact that 
$\omega(c)$ maps to $c(F^*\otimes L - E)$, proves Proposition~\ref{tilde}.

Given the setup of the degeneracy loci section, the perpendiculars of the bundles define 
new flags  $G_\bull$ in $V$ and $H_\bull$ in $W$, with $G_0$ and $H_0$ of the same rank;  namely, set  $G_i = F_{-i}^\perp = \Hom(W/F_{-i},L)$ and $H_i = E_{-i}^\perp = \Hom(V/E_{-i},L)$.  From the formula 
\[ 
 k_{\omega(w)}(p,q) \, = \, k_{w^{-1}}(-q,-p) ,
\]
and the fact that the kernels of  $F_q \to \Hom(E_p,L)$  and $H_{-p} \to G_{-q}$ are both $F_q \cap H_{-p}$, it follows that 
\[
 \Omega_w(E_\bull, F_\bull) \, = \,   \Omega_{\omega(w^{-1})}(G_\bull, H_\bull) .
\]
Proposition~\ref{omega} then follows from the duality formula in Proposition~\ref{duality}.

Consider the case where $E = V$ and $F = W$, with the bilinear form non-degenerate.  Set $G_i = F_i^\perp$, and $H_i = E_i^\perp$.  Note that $G_i$ and $H_i$ have rank $i$, and $c_1(G_i/G_{i-1}) = z-y_i$, and   
$c_1(\Ker( \Hom(H_i,L) \to \Hom(H_{i-1},L) ) = x_i$.  We have 
\[ 
G_1  \hookrightarrow \cdots \hookrightarrow G_N = E  \stackrel{\cong}{\rightarrow} \Hom(F,L) = \Hom(H_N,L)  \twoheadrightarrow  \cdots  \twoheadrightarrow \Hom(H_1,L) .
\]
We have the locus $\Omega_w^{\text{old}}$, given by requiring the rank of $G_q \to \Hom(H_p,L)$ to be at most 
$r_w(p,q) =  \#\{a  \leq p \mid w(a) \leq q \}$ for all $p$ and $q$.  By \cite{F1}, the formula for this locus is 
\[
 [\Omega_w^{\text{old}}] = \frak{S}_w(x,z-y).
\]
In fact, this locus $\Omega_w^{\text{old}}$ is the same as the locus denoted above by 
$\Omega_w(E_\bull,F_\bull)$.  (This is an easy consequence of the ``rank-nullity theorem":
\[
\dim(E_p \cap F_q^\perp) + p \; = \; \dim(E_p^\perp \cap F_q) + q, 
\]
and the fact that $r_w(p,q) + k_w(p,q) = p$.)
The nondegeneracy condition implies that $c = c(E^*\otimes L - F) = 1$, so 
\[
 [\Omega(E_\bull,F_\bull)] = \tenS_w(0,x,y) .
\]
This gives a geometric proof of the fact that 
\[
\tenS_w(0,x,y) = \frak{S}_w(x,z-y).
\]
The same idea gives a geometric proof of Theorem~\ref{DefThm}.

The back-stable identity for twisted enriched Schubert polynomials can also be motivated and proved from degeneracy loci.  Given filtrations $E_\bull$ and $F_\bull$  in $V$ and $W$, and any integer $m$, there are filtrations $\gamma^m(E_\bull)$ and $\gamma^m(F_\bull)$, defined by setting $\gamma^m(E_i) = E_{m+i}$ and $\gamma^m(F_i) = F_{m+i}$.  Then 
\[
 \Omega_{\gamma^m(w)}(\gamma^m(E_\bull),\gamma^m(F_\bull)) = \Omega_w(E_\bull,F_\bull).
\]
For positive $m$, 
\begin{align*}
 c(E_m^*\otimes L - F_m ) = c( E_0^*\otimes L - F_0)\cdot c(F_0/F_m)/c((E_0/E_m)^*\otimes L),
\end{align*}
which agrees with the formula  $\gamma^m(c) = c \cdot \prod_{i=1}^m \frac{1+y_i}{1+z-x_i}$, so the back-stable formula of Proposition~\ref{p.back-stable} follows.  (Similarly for 
negative $m$.)

To see how the ``intersection loci'' of Proposition~\ref{altdegloci} are related to the ``tensor rank'' loci defined in Theorem~\ref{degloci}, given $G_\bull$ and $F_\bull$ in $W$ as in Proposition~\ref{altdegloci}, set $V = W^*$, 
so we have the canonical (non-degenerate) bilinear form from $V \times W$ to $L$, where $L$ is the trivial line bundle.  Set 
$E_p = (W/G_p)^* \subset V$.  Then $\Ker(F_q \to E_p^*) = G_p \cap F_q$, so the loci are equal: $\Omega_w(G_\bull \cap F_\bull) = \Omega_w(E_\bull, F_\bull)$.  Since $E^* = V/G$, $c(V - G - F) = c(E^* - F)$.  And $G_i/G_{i-1} = (E_{i-1}/E_i)^*$, so $c_1(G_i/G_{i-1}) = -x_i$.  One can also formulate this intersection loci formula with a line bundle, replacing duals by $L$-duals, in which case $c_1(G_i/G_{i-1}) = z-x_i$.

For Corollary~\ref{Pawlowski},
let $\mathbb{W} = W \oplus W = U \oplus U \oplus V \oplus V$, with flags defined by
\[
  \mathbb{G}_p = U\oplus S \oplus G_p \quad \text{and} \quad \mathbb{F}_q = \Delta( F_q\oplus V ),
\]
where $\Delta\colon W \hookrightarrow \mathbb{W}$ is the diagonal.  Then
\[
  \mathbb{G}_p \cap \mathbb{F}_q = \Delta( S \cap F_q \oplus G_p ),
\]
so $\Gamma\Omega_w = \Omega_w(\mathbb{G}_\bull \cap \mathbb{F}_\bull)$ and the formula follows.

Our proof of Theorem~\ref{AtoC} is geometric, although one could also deduce it from localization.   Start with a vector bundle $V$ of large rank $2m$ on a (nonsingular) variety $Y$, with 
an an $L$-valued symplectic form $V \times V \to L$.  To be assured that the situation is general enough, take $X$ to consist of
pairs of complete flags in $V$; that is, 
take  $X \to Y$ to be the fiber product of two copies of the complete flag bundle of isotropic flags of $V$: $X = \Fl(V) \times_Y \Fl(V)$.  
Take  $E_\bull$ and $F_\bull$ to be the pullbacks of the tautological flags on the two factors (writing as usual $V$ for the pullback of $V$ to $X$, and index these flags so $E_i$ and $F_i$ have ranks $m - i$ for $-m \leq i \leq m$.
There is a subvariety $X_C$ of $X$, consisting of those pairs of flags such 
that $E = E_0$ and $F = F_0$ are isotropic, and, in addition, 
\[ E_i = (E_{-i})^\perp   \;\; \text{ and } \;\;  F_i = (F_{-i})^\perp \]
for all $i$.  To match the type C notation, set, for $i \geq 1$,
\[ \bar{E}_i = E_{i-1}   \;\; \text{ and } \;\; \bar{F}_i = F_{i-1}, \]
and $\bar{E}_i = E_i$ and $\bar{F_i} = F_i$ for $i < 0$.
  
Using the notation $\Omega_w^A(E_\bull,F_\bull)$ and $\Omega_w^C(\overline{E}_\bull,\overline{F}_\bull)$ for the two degeneracy loci, the theorem follows from the
\begin{claim*} For 
any $w$ in $S_m \subset S_+$, the intersection of $\Omega_w^A(E_\bull,F_\bull)$ with $X_C$ is 
$\Omega_w^C(\bar{E}_\bull,\bar{F}_\bull)$.
\end{claim*}

For $w$ in $S_m$, the type C locus $\Omega_w^C(\bar{E}_\bull,\bar{F}_\bull)$ is given by the conditions
 \[ \dim(\bar{E}_{p+1} \cap \bar{F}_{-q}) \geq  \# \{ a \geq p+1 \mid - w(a) \geq - q \} \]
 for $p \geq 0$ and $q \geq 1$.  Since $\bar{F}_{-q} = F_q^\perp$, this says that the kernel of nullity of  $F_q \to \Hom(E_p, L)$ has dimension at least $ \# \{ a > p \mid w(a) \leq  q \}$, which is exactly the type A condition to be in $\Omega_{w^{-1}}^A(F_\bull,E_\bull) = \Omega_w^A(E_\bull,F_\bull)$. 
 
\section{Table of Type A Enriched Schubert Polynomials}

The following gives a table of the type A Schubert polynomials for  permutations that are the identity outside the interval $[0,3]$.  Here  $S_\lambda$ denotes $S_\lambda(c)$; for the twisted polynomials  $\tenS_w$, replace each $x_i$ by $x_i - z$. 

\vspace{.2cm}

{\small
\begin{tabular}{|l|l|}
\hline

$w$ & $\enS_w(c,x,y)$ \\ \hline \hline

$0\;1\;2\;3$ & $ 1 $ \\ \hline
$0\;1\;3\;2$ & $ x_1+x_2+y_1+y_2+S_{{1}} $ \\ \hline
$0\;2\;1\;3$ & $ x_1+y_1+S_{{1}} $ \\ \hline
$0\;2\;3\;1$ & $ (x_1+y_1)(x_2+y_1) + \left( x_1+ x_2+ y_1 \right) S_{{1}}+S_{{1,1}} $ \\ \hline
$0\;3\;1\;2$ & $  (x_1+y_1)(x_1+y_2) + \left(  x_1+ y_1+ y_2 \right) S_{{1}}+S_{{2}} $ \\ \hline
$0\;3\;2\;1$ & $\begin{array}{l}  (x_1+y_1)(x_1+y_2)(x_2+y_1)  + \left(  x_1+ y_1+y_2 \right)  \left( x_1+ x_2+y_1 \right) S_{{1}} \\ \quad + \left( x_1+ x_2+ y_1 \right) S_{{2}}+ \left(  x_1+ y_1+ y_2 \right) S_{{1,1}}+S_{{2,1}}  \end{array}$ \\ \hline
$1\;0\;2\;3$ & $ S_{{1}} $ \\ \hline
$1\;0\;3\;2$ & $  \left(  x_1+ x_2+y_1+y_2 \right) S_{{1}}+S_{{2}}+S_{{1,1}} $ \\ \hline
$1\;2\;0\;3$ & $ x_1\,S_{{1}}+S_{{1,1}} $ \\ \hline
$1\;2\;3\;0$ & $ x_1\,x_2\,S_{{1}}+ \left( x_1+x_2
 \right) S_{{1,1}}+S_{{1,1,1}} $ \\ \hline
$1\;3\;0\;2$ & $ x_1 \left( x_1+y_1+y_2 \right) S_{{1}} + x_1\,S_{{2}} + \left( x_1+y_1+y_2 \right) S_{{1,1}}+S_{{2,1}} $ \\ \hline
$1\;3\;2\;0$ & $ \begin{array}{l} x_1 x_2 \left( x_1+y_1+y_2 \right) \,S_{{1}}+x_1\,x_2\,S_{{2}} + \left( x_1+y_1+y_2 \right) \left( x_1+x_2 \right) S_{{1,1}} \\ \quad + \left( x_1+x_2 \right) S_{{2,1}} + \left( x_1+y_1+y_2 \right) S_{{1,1,1}}+S_{{2,1,1}} \end{array} $ \\ \hline
$2\;0\;1\;3$ & $ y_1\,S_{{1}}+S_{{2}} $ \\ \hline
$2\;0\;3\;1$ & $  \left( x_1+x_2+y_1 \right) y_1\,S_{{1}} + \left( x_1+x_2+y_1 \right) S_{{2}}+y_1\,S_{{1,1}} +S_{{2,1}} $ \\ \hline
$2\;1\;0\;3$ & $ x_1\,y_1\,S_{{1}}+ x_1 S_{{2}} +y_1\,S_{{1,1}}+S_{{2,1}} $ \\ \hline
$2\;1\;3\;0$ & $ \begin{array}{l} x_1\,x_2\,y_1\,S_{{1}} + x_1\,x_2\,S_{{2}} + \left( x_1+x_2 \right) y_1\,S_{{1,1}} +\left( x_1+x_2 \right) S_{{2,1}} \\ \quad  +y_1\,S_{{1,1,1}} +S_{{2,1,1}} \end{array} $ \\ \hline
$2\;3\;0\;1$ & $ \begin{array}{l} \left( x_1+y_1 \right) x_1\,y_1\,S_{{1}} + \left( x_1+y_1 \right) x_1\,S_{{2}} + \left( x_1+y_1 \right) y_1\,S_{{1,1}} \\ \quad + \left( x_1+y_1 \right) S_{{2,1}}+S_{{2,2}} \end{array} $ \\ \hline
$2\;3\;1\;0$ & $ \begin{array}{l} \left( x_1+y_1 \right) x_1\,x_2\,y_1\,S_{{1}}  + \left( x_1+y_1 \right) x_1\,x_2\,S_{{2}} + \left( x_1+y_1 \right)  \left( x_1+x_2 \right) y_1\,S_{{1,1}} \\ \quad + \left( x_1+y_1 \right)  \left( x_1+x_2 \right) S_{{2,1}}  + \left( x_1+y_1 \right) y_1\,S_{{1,1,1}}  \\ \quad + \left( x_1+x_2 \right) S_{{2,2}} + \left( x_1+y_1 \right) S_{{2,1,1}}+S_{{2,2,1}} \end{array} $ \\ \hline
$3\;0\;1\;2$ & $ y_1\,y_2\,S_{{1}}+ \left( y_1+y_2 \right) S_{{2}}+S_{{3}} $ \\ \hline
$3\;0\;2\;1$ & $  \begin{array}{l} \left( x_1+x_2+y_1 \right) y_1\,y_2\,S_{{1}} + \left( y_1+y_2 \right)  \left( x_1+x_2+y_1 \right) S_{{2}} +y_1\,y_2\,S_{{1,1}} \\ \quad + \left( x_1+x_2+y_1 \right) S_{{3}}+ \left( y_1+y_2 \right) S_{{2,1}}+S_{{3,1}} \end{array} $ \\ \hline
$3\;1\;0\;2$ & $ \begin{array}{l} x_1\,y_1\,y_2\,S_{{1}}+ \left( y_1+y_2 \right) x_1\,S_{{2}} +y_1\,y_2\,S_{{1,1}}+x_1\,S_{{3}} \\ \quad +\left( y_1+y_2 \right) S_{{2,1}}+S_{{3,1}} \end{array} $ \\ \hline
$3\;1\;2\;0$ & $ \begin{array}{l} x_1\,x_2\,y_1\,y_2\,S_{{1}} + \left( y_1+y_2 \right) x_2\,x_1\,S_{{2}} + \left( x_1+x_2 \right) y_1\,y_2\,S_{{1,1}} \\ \quad +x_1\,x_2\,S_{{3}}+ \left( y_1+y_2 \right)  \left( x_1+x_2 \right) S_{{2,1}} +y_1\,y_2\,S_{{1,1,1}}  \\ \quad + \left( x_1+x_2 \right) S_{{3,1}}+ \left( y_1+y_2 \right) S_{{2,1,1}}+S_{{3,1,1}} \end{array} $ \\ \hline
$3\;2\;0\;1$ & $ \begin{array}{l} \left( x_1+y_1 \right) x_1\,y_1\,y_2\,S_{{1}} + x_1 \left( x_1+y_1 \right)  \left( y_1+y_2 \right) S_{{2}}  + \left( x_1+y_1 \right) y_1\,y_2\,S_{{1,1}}\\ \quad + \left( x_1+y_1 \right) x_1\,S_{{3}}+ \left( y_1+y_2\right)  \left( x_1+y_1 \right) S_{{2,1}} + \left( x_1+y_1 \right) S_{{3,1}} \\ \quad  + \left( y_1+y_2 \right) S_{{2,2}} + S_{{3,2}}  \end{array} $ \\ \hline
$3\;2\;1\;0$ & $ \begin{array}{l} \left( x_1+y_1 \right) x_1\,x_2\,y_1\,y_2\,S_{{1}} + \left( y_1+y_2 \right)  \left( x_1+y_1 \right) x_1\,x_2\,S_{{2}} \\ \quad + \left( x_1+y_1 \right)  \left( x_1+x_2 \right) y_1\,y_2\,S_{{1,1}} + \left( x_1+y_1 \right) x_1\,x_2\,S_{{3}}   \\ \quad  + \left( y_1+y_2 \right)  \left( x_1+y_1 \right) 
 \left( x_1+x_2 \right) S_{{2,1}} + \left( x_1+y_1 \right) y_1\,y_2\,S_{{1,1,1}} \\ \quad + \left( x_1+y_1 \right)  \left( x_1+x_2 \right) S_{{3,1}} + \left( y_1+y_2 \right)  \left( x_1+x_2 \right) S_{{2,2}} \\ \quad+ \left( y_1+y_2 \right)  \left( x_1+y_1 \right) S_{{2,1,1}} + \left( x_1+x_2 \right) S_{{3,2}} + \left( x_1+y_1 \right) S_{{3,1,1}} \\ \quad + \left( y_1+y_2 \right) S_{{2,2,1}}  +S_{{3,2,1}}  \end{array} $ \\ \hline

\end{tabular}

}


\end{document}